\documentclass[11pt]{article}
\usepackage{latexsym}
\usepackage{amssymb,latexsym}
\usepackage{amsmath}
\usepackage{mathrsfs}
\usepackage{amscd}
\usepackage{amsfonts}
\usepackage[notref,notcite]{showkeys}
\newcommand{\noin}{\noindent}

\newcommand{\grad}{{{\rm{grad}\,}}}

\newcommand{\Hess}{{{\rm Hess}\,}}

\newtheorem{theorem}{\bf Theorem}[section]

\newtheorem{proposition}[theorem]{\bf Proposition}

\begin{document}

\title{Complete submanifolds of $\mathbb{R}^{n}$ with finite topology}

\author{G. Pacelli Bessa  \and Luqu\'esio Jorge \and
 J. F\'{a}bio Montenegro}
\date{}
\maketitle
\begin{abstract}
We show that a complete $m$-dimensional immersed submanifold $M$ of
$\mathbb{R}^{n}$  with $a(M)<1$ is properly immersed and have finite
topology, where $a(M)\in [0,\infty]$ is an scaling invariant number
that gives the rate  that the norm of the second fundamental form
decays to zero at infinity. The class of submanifolds $M$ with
$a(M)<1$ contains all complete minimal surfaces in $\mathbb{R}^{n}$
with finite total curvature, all $m$-dimensional minimal
submanifolds
 $M $ of $ \mathbb{R}^{n}$ with finite
total scalar curvature $\smallint_{M}\vert \alpha \vert^{m}
dV<\infty $ and all complete $2$-dimensional complete surfaces with
$\smallint_{M}\vert \alpha \vert^{2} dV<\infty $ and nonpositive
curvature with respect to every normal direction, since $a(M)=0$ for
them.

\end{abstract}
\section{Introduction}
 \hspace{2mm} Let $M$   be a  complete
  surface minimally immersed in $\mathbb{R}^{n}$ and let
 $K$
  be
 its Gaussian  curvature.  Osserman  \cite{osserman} for $n=3$  and Chern-Osserman
 \cite{chern-osserman} for $n\geq 3$
   proved that
$\smallint_{M}\vert K\vert dV<\infty $ if and only if $M$ is
conformally equivalent to a compact Riemann surface $\overline{M}$
punctured at a finite number of points $\{p_{1},\ldots,p_{r}\}$ and
the Gauss map $\Phi:M\to \mathbb{G}_{2,n}$ extends to a holomorphic
map $\overline{\Phi}:\overline{M}\to \mathbb{G}_{2,n}$, see
\cite{lawson} for a beautiful exposition. M.  Anderson
\cite{anderson} proved a higher dimension version of Chern-Osserman
finite total curvature theorem, i.e. a complete $m$-dimensional
minimally immersed submanifold $M$ of $\mathbb{R}^{n}$ has finite
total scalar curvature $\smallint_{M}\vert \alpha \vert^{m}
dV<\infty $ if and only if $M$ is $C^{\infty}$- diffeomorphic to a
compact smooth Riemannian manifold $\overline{M}$ punctured at a
finite number of points $ \{p_{1},\ldots,p_{r}\}$
 and the Gauss map $\Phi$ on $M$
extends to a $C^{\infty}$- map $\overline{\Phi}$ on $\overline{M}$,
where $\vert \alpha \vert$ is the norm of the second fundamental
form of $M$.

\hspace{-3mm} It is interesting that   these results have
appropriate versions in the non-minimal setting. B.  White
\cite{white}, proved that  a complete $2$-dimensional
 surface $M$ immersed in $\mathbb{R}^{n}$ with    $\smallint_{M}\vert
\alpha \vert^{2} dV<\infty $ and  nonpositive  curvature with
respect to every normal direction\footnote{$M\subset \mathbb{R}^{n}$
is nonpositively curved with respect to each normal direction  at
$x$ if ${\rm det} (\eta \cdot \alpha ( ,))\leq 0$ for all normals
$\eta$ to $M$ at $x$, see \cite{white}.} is homeomorphic to a
compact Riemann surface $\overline{M}$ punctured at finite number of
points  $\{p_{1},\ldots,p_{r}\}$, its Gauss map $\Phi$ extends
continuously to all of $\overline{M}$ and $M$ is properly immersed.
It should be observed that the  properness of $M$ in White's Theorem
is consequence of from the first two statements about the immersion,
i.e.  Jorge and  Meeks \cite{jorge-meeks},  proved that  a complete
$m$-dimensional
 immersed submanifold $M$ of $\mathbb{R}^{n}$, homeomorphic to a
compact Riemann manifold $\overline{M}$ punctured at finite number
of points $\{p_{1},\ldots,p_{r}\}$ and such that the Gauss map
$\Phi$ extends continuously to all of $\overline{M}$ is properly
immersed.

S. Muller and V. Sverak \cite{muller-sverak}, answering a question
of White, proved that a complete $2$-dimensional
 surface $M$ immersed in $\mathbb{R}^{n}$ with    $\smallint_{M}\vert
\alpha \vert^{2} dV<\infty $ is properly immersed.
 For higher
dimension submanifolds,  Do Carmo and  Elbert \cite{doCarmoElbert}
proved that a
 complete hypersurface  $M$ of $ \mathbb{R}^{n+1}$ with strong finite total scalar curvature, ($\smallint_{M}(\vert
\alpha \vert ^{n-1}+\vert \nabla \alpha\vert^{n-2})dV<\infty$) is
proper and it is $C^{\infty}$- diffeomorphic to a compact smooth
manifold $\overline{M}$ punctured at a  finite number of points
$\{p_{1},\ldots,p_{r}\}$. Moreover, they showed that if
 the Gauss-Kronecker curvature does not
change sign in a punctured neighborhood of each $p_{i}$,
$i=1,\ldots,r$ then the Gauss map $\Phi $ on $M$ extends
continuously to all of $\overline{M}$.

\hspace{-3mm} The purpose of this paper is to put another piece on
this puzzle showing  that a complete  $m$-dimensional submanifold of
$\mathbb{R}^{n}$ with
 the norm of the second fundamental form uniformly  decaying to zero $\vert \alpha(x)\vert \to 0$
  as $x\to \infty$  in a certain rate is proper
and has finite topological type. The decaying rate of  $\vert \alpha
(x)\vert \to 0$   considered is not fast enough to make
$\smallint_{M}\vert \alpha (x)\vert^{m}dV <\infty$.

 To be more precise, let
 $M$ be a complete $m$-dimensional
submanifold of $\mathbb{R}^{n}$  and let $K_{1}\subset K_{2}\subset
\ldots$ be an exhaustion sequence of $M$ by compact sets. Fix a
point
  $p\in K_{1}$ and set
$a_{i}=\sup\{\rho (x)\cdot \vert \alpha(x)\vert, \,x\in M\setminus
K_{i} \}$, where $\rho (x)={\rm dist}_{M}(p,x)$ and $\vert\alpha (x)
\vert$ is the norm of the second fundamental form of $M$ at $x$. The
$a_{i}$'s form a non-increasing sequence $\infty\geq a_{1}\geq
a_{2}\geq \cdots\geq 0$ with $a_{1}=\infty$ iff $a_{l}=\infty$ for
all $l\geq 1$.  Define the (possibly extended) scaling invariant
number $a(M)=\lim_{i\to \infty}a_{i}\in [0,\infty]$. It can be shown
that $a(M)$ does not depend on the exhaustion sequence nor on the
point $p$. It follows from the work of Jorge-Meeks
\cite{jorge-meeks} that complete $m$-dimensional Riemannian
submanifolds $M$ of $ \mathbb{R}^{n}$homeomorphic to a compact
Riemann manifold $\overline{M}$ punctured at finite number of points
$\{p_{1},\ldots,p_{r}\}$  and well defined normal vector at infinity
have $a(M)=0$. In particular,  complete minimal surfaces in
$\mathbb{R}^{n}$ with finite total curvature, complete
$2$-dimensional complete surfaces with $\smallint_{M}\vert \alpha
\vert^{2} dV<\infty $ and nonpositive curvature with respect to
every normal direction considered by White or the $m$-dimensional
minimal submanifolds
 $M $ of $ \mathbb{R}^{n}$ with finite
total scalar curvature considered by Anderson have $a(M)=0$. In our
main result, we prove that the larger class   of complete
$m$-dimensional immersed submanifolds of $\mathbb{R}^{n}$ with
$a(M)<1$ share some properties with theses submanifolds with
$a(M)=0$. We prove the following theorem.

\vspace{3mm}
\begin{theorem}\label{theorem1} Let $M$
 be a complete $m$-dimensional   submanifold of $\mathbb{R}^{n}$ with
$a(M)<1$. Then $M$ is properly immersed and it is $C^{\infty}$-
diffeomorphic to a  compact smooth manifold $\overline{M}$ with
boundary.
\end{theorem}
 Observe that
$\smallint_{M}\vert \alpha \vert^{m} dV<\infty $ is not equivalent
to $a(M)<1$. However, one might ask if Theorem (\ref{theorem1})
holds under finite total scalar curvature $\smallint_{M}\vert \alpha
\vert^{m} dV<\infty $.

\vspace{3mm}

\hspace{-3mm} For  complete $m$-dimensional minimal submanifolds $M$
of $\mathbb{R}^{n}$ we define the increasing sequence
$b_{i}=\inf\{\rho^{2} (x)\cdot Ric (x)(\nu,\nu), \, \vert \nu \vert
=1, \,x\in M\setminus K_{i} \}$ with $b_{1}=-\infty $ iff
$b_{l}=-\infty $ for all $l\geq 1$. Define the scaling invariant
number $b(M)=\lim_{i\to \infty}b_{i}\in [-\infty, 0]$. Again, it can
be shown that $b(M)$ does not depend on the exhaustion sequence nor
on the point $p$. The proof of  Theorem (\ref{theorem1}) can be
slightly  modified to prove the following version for minimal
submanifolds.
\begin{theorem}\label{theorem2} Let $M$
 be a complete $m$-dimensional  minimal submanifold of $\mathbb{R}^{n}$
 with $b(M)>-1$. Then  $M $ is properly immersed and it is
$C^{\infty}$- diffeomorphic to a  compact smooth manifold
$\overline{M}$ with boundary.
\end{theorem}

\section{Proof of Theorem \ref{theorem1}}
\subsection{$M$ is properly immersed}
 Let $\varphi :M^{m} \hookrightarrow \mathbb{R}^{n}$ be a
complete  submanifold  with $a(M)<1$ and let $p\in M$ be a fixed
point such that $\varphi(p)=0\in \mathbb{R}^{n}$. There exist a
geodesic  ball $B_{M}(p,R_{0})$ centered at $p$ with radius $R_{0}$
such that for all $x\in M\setminus B_{M}(p,R_{0})$ we have that
$\rho (x)\vert \alpha(x)\vert\leq c<1$. Let $f:M^{m}\to \mathbb{R}$
given by $f(x)=\vert \varphi (x)\vert^{2}$. Fix a point $x\in
M\setminus B_{M}(p,R_{0})$ then for $\nu\in T_{x}M$, $\vert \nu
\vert =1$ we have that
\begin{eqnarray}\label{eqHess} \frac{1}{2}\,\Hess f (x)(\nu,\nu)& =&
1
+\langle \varphi(x),\,\alpha (x) (\nu,\nu)\rangle \nonumber \\
&\geq & 1-\vert \varphi (x)\vert \cdot \vert \alpha (x)
 \vert \\
 &\geq & 1-\rho (x)\vert \alpha (x)\vert \nonumber \\
 &\geq & 1-c\nonumber
\end{eqnarray}
 Let $\sigma :[0,\rho (x)]\to M^{m}$ be a minimal
geodesic from $p$ to $x$. We have from (\ref{eqHess}) that for all
$t\geq R_{0}$ that
 $(f\circ\sigma)''(t))=\Hess f (\sigma
(t))(\sigma', \sigma')\geq 2(1-c)$ and for $t<R_{0}$ that
$(f\circ\sigma)''(t))\geq b$, $b=\inf_{x\in B_{M}(p,R_{0})} \{ \Hess
f (x)(\nu,\nu), \,\vert \nu\vert=1\}$. Thus
\begin{eqnarray}\label{eq2}
(f\circ\sigma) ' (s)&=& \int_{0}^{s}(f\circ \sigma)'' (\tau)d\tau \nonumber  \\
&\geq & \int_{0}^{R_{0}}b\,d\tau + \int_{R_{0}}^{s}(1-c)d\tau\\
&&\nonumber \\ &\geq & b\,R_{0}+ (1-c) (s-R_{0})\nonumber
\end{eqnarray}
\begin{eqnarray}
 f(x)&=& \int_{0}^{\rho(x)}(f\circ \sigma)' (s) \, ds\nonumber \\
 &\geq & \int_{0}^{\rho(x)}b\,R_{0}+ (1-c) (s-R_{0})\,ds \label{eq3} \\
 &=&b\,R_{0}\,\rho (x) + (1-c)\left( \frac{\rho(x)^{2}}{2}-R_{0}\rho(x)\right)\nonumber
\end{eqnarray} Thus $\varphi (x)\vert^{2}\geq
(b-1+c)R_{0}\rho(x)+(1-c)\rho(x)^{2}/2$ for all $x\in M\setminus
B_{M}(p,R)$. In fact, this proof that $f$ is proper proves that
following proposition.
\begin{proposition}\label{propositionProper} Let  $f:M\to
\mathbb{R}$ be  a $C^{2}$-function defined on a complete Riemannian
manifold such that $\Hess f(x) \geq g(\rho (x))$, where $\rho $ is
the distance function  to $x_{o}$ and $g:[0,\infty)\to \mathbb{R}$
is a piecewise continuous function. Setting $G(t)=f(x_{o})-\vert
\grad f (x_{o})\vert \,t +
\smallint_{0}^{t}\smallint_{0}^{s}g(u)duds$, $t\in [0,\infty)$ we
have that if $G$ is proper and bounded from below
  then $f$ is proper.
\end{proposition}
\subsection{$M$ has finite topology}
Let $\varphi :M^{m} \hookrightarrow \mathbb{R}^{n}$ be a complete
properly immersed submanifold  with $a(M)<1$. Let $p\in M$ be  such
that $\varphi(p)=0\in \mathbb{R}^{n}$ and suppose without loss of
generality that $R(x)=\vert \varphi (x)\vert$, $x\in M$ is a Morse
function. Let $r>0$ be such that $\Gamma_{r}=\varphi (M)\,\cap\,
\mathbb{S}^{n-1}(r)$ is a compact submanifold of
$\mathbb{S}^{n-1}(r)$ and $\rho (x)\cdot \vert \alpha (x)\vert \leq
c<1$ for all $x\in M\setminus \varphi^{-1}(B_{\mathbb{R}^{n}}(r))$.
For $y\in \Gamma_{r}$, let $\eta (y)=y/r$, the unit vector
perpendicular to $T_{y}\mathbb{S}^{n-1}(r)$. There is only one unit
vector
 $\nu (y)\in T_{y}\varphi (M) $  such that
 $T_{y}\varphi(M)=T_{y}\Gamma_{r}\oplus [[\nu (y)]]$
 and $\langle \nu (y), \eta (y)\rangle > 0$. Here $T_{y}\varphi(M)$,
 $T_{y}\mathbb{S}^{n-1}(r)$ are the tangent spaces of $\varphi (M)$
 and $\mathbb{S}^{n-1}(r)$ at $y$ respectively and
 $[[\nu (y)]] $ is the vector
  space generated by $ \nu (y)$. This procedure defines
 a smooth vector field  $\nu (y)$ on the neighborhood $\varphi^{-1}(B_{\mathbb{R}^{n}}(r+\epsilon)\setminus
 B_{\mathbb{R}^{n}}(r-\delta))$ of $\varphi^{-1}(\Gamma_{r})$  for
 some
 $\epsilon, \delta
 >0$.
  Consider the function $\psi$
  defined on that neighborhood given by
$$\psi (y)=\langle \nu(y), \eta (y)\rangle=\cos \theta (y)$$

\noin Identify $X\in TM$  with $d\varphi X$. We have that
$X(R)=\langle X,\eta\rangle$ and  writing $\eta(y)= \sin
\theta(y)\,\nu^{\ast} (y) + \cos \theta (y)\,\nu (y)$, $ \nu^{\ast}
(y)\perp \nu (y)$, we have that $ \grad R=\psi\,\nu$. Now for each
$y\in \Gamma_{r}$ consider $\xi(t,y)$ the solution of the following
problem on $\varphi (M)$
\begin{equation}\label{eqData}
\left\{ \begin{array}{l} \xi_t=\displaystyle\frac{1}{\psi}\,\nu(\xi) \\ \\
\xi(0,y)=y
\end{array}\right.
\end{equation}
 We have   for $R(t,y)=\vert\xi(t,y)\vert$ that
$$\displaystyle
R_t=\langle \grad R,\frac{1}{\psi}\,\nu \rangle = \langle \psi \nu
,\frac{1}{\psi} \,\nu \rangle =1 \Leftrightarrow R=R(t,y)=t+r
$$
\noin  We will  derive a differential equation that the function
$\psi \circ \xi (t,y)$ satisfies.
$$\displaystyle
\psi_t=\xi_t\langle \nu ,\eta \rangle = \langle D_{\frac{1}{\psi}\nu
}\nu ,\eta \rangle + \langle \nu , D_{\xi_t}\eta \rangle = \langle
\nabla_{\nu}\nu +\alpha(\nu ,\nu ),\eta \rangle +\langle \nu,
D_{\xi_t}(\frac{\xi}{R}) \rangle
$$
But $\langle \nu ,\nu \rangle =1 \; \Rightarrow \langle \nu
,\nabla_{\nu}\nu \rangle =0$ and $\nabla_{\nu}\nu \in T_xM
\;\Rightarrow$ $\nabla_{\nu}\nu \in (T_xM \cap T_xS_R^n)$
$\Rightarrow \langle \nabla_{\nu}\nu ,\eta\rangle =0.$ On other
hand
$$
D_{\xi_t}(\frac{\xi}{R})=\frac{\frac{1}{\psi}\nu}{R} -
\frac{R_t}{R^2}\varphi = \frac{1}{R\psi}\nu -\frac{1}{R}\eta
$$
then
\begin{equation}\label{eqdf}
\psi_t=\frac{1}{\psi}\langle \alpha(\nu ,\nu ),\eta \rangle
+\frac{1}{\psi
R}-\frac{\psi}{R}=\frac{\sqrt{1-\psi^2}}{\psi}\langle \alpha (\nu,
\nu ), \nu^* \rangle +\frac{1-\psi^2}{\psi R}
\end{equation}
To  determine a differential equation satisfied by $\sin
\theta(t,y)=\sqrt{1-\psi^2}$ we proceed as follows. By (\ref{eqdf})
we have

\begin{equation}\label{eqdf2}
\frac{\psi \psi_t}{\sqrt{1-\psi^2}}=\langle \alpha (\nu, \nu ),
\nu^* \rangle +\frac{\sqrt{1-\psi^2}}{R}
\end{equation}
Observing that $R(t,y)=t+r$, equation  (\ref{eqdf2}) can be written
as
\begin{equation}-(t+r)(\sqrt{1-\psi^2})_{t}=(t+r)\langle \alpha (\nu, \nu ),
\nu^* \rangle +\sqrt{1-\psi^2}\label{eqdf3}
 \end{equation} and rewritten  as
 \begin{equation}\left[(t+r)\sqrt{1-\psi^2})\right]_{t}+(t+r)\langle \alpha (\nu, \nu ),
\nu^* \rangle=0\label{eqdf4}
 \end{equation}
Integrating the equation (\ref{eqdf4}) we have the following
equation
\begin{equation}\sqrt{1-\psi^{2}}=\frac{r}{t+r}\sqrt{1-\psi_{o}^{2}}
-\frac{1}{t+r} \int_{0}^{t}(s+r)\langle \alpha (\nu, \nu ), \nu^*
\rangle\,ds\label{eqdf5}
\end{equation}
where $\psi_{o}=\psi(\xi(0,y))$. Since $ \sin \theta (\xi(t,y))=
\sqrt{1-\psi^{2}}$ we rewrite (\ref{eqdf5}) in the following form
\begin{equation}\sin \theta
(\xi(t,y))=\frac{r}{t+r}\sin\theta (\xi (0,y))-\frac{1}{t+r}
\int_{0}^{t}(s+r)\langle \alpha (\nu, \nu ), \nu^*
\rangle\,ds\label{eqdf6}
\end{equation}
  Now, $$-\langle \alpha
(\nu, \nu),\nu^{\ast}\rangle (\xi(s,y)) \leq \vert \alpha \vert (\xi
(s,y)) \leq c/\rho (\xi(s,y))\leq c/ R(s,y).$$ Substituting in
(\ref{eqdf6}) and noticing that $ R(s,y)=s+r$ we have that
\begin{eqnarray}\sin \theta
(\xi(t,y))& \leq & \frac{r}{t+r}\sin\theta (\xi (0,y))+\frac{1}{t+r}
\int_{0}^{t}(s+r)\frac{c}{s+r}\,ds\label{eqdf7}\nonumber
\\
&& \nonumber \\
&=& \frac{c\,t+r\sin\theta (\xi (0,y))}{t+r}<1,\,\, \forall t\geq 0.
\end{eqnarray}
 Recall that
the critical points of $R$ are those $x$ such that $\psi (x)=0$, or
those points where $\sin \theta(x)=1$. Thus, along the integral
curves $\xi (t,y)$, $y\in \Gamma_{r}$ there is no critical point for
the function $R(x)=\vert \varphi (x)\vert$. This shows that outside
the compact set $M\setminus B_{M}(p,r)$ there are no critical points
for $R$. Since $R$ is a Morse function, its  critical points are
isolated  thus there are finitely many of them. This shows that $M$
has finite topology.

\section{Sketch of proof for Theorem \ref{theorem2}}
 Let $\varphi :M^{m} \hookrightarrow
\mathbb{R}^{n}$ be a complete minimal submanifold and let $x\in M$,
$\nu \in T_{x}M$ and $\{e_{1},\ldots e_{m}=\nu\}$ an orthonormal
basis for $T_{x}M$. Using the Gauss equation we can compute the
Ricci curvature in the direction $\nu$ by
\begin{eqnarray}\label{eqRic}
Ric(x) (\nu )& =& \langle \sum_{i=1}^{m} \alpha_{ii}
,\alpha_{mm}\rangle-\sum_{i=1}^{m-1}\vert\alpha_{im}\vert^2\nonumber
\\&  =& \langle mH-\alpha_{mm},\alpha_{mm}\rangle
-\sum_{i=1}^{m-1}|\alpha_{im}|^2
\\ & =&-\sum_{i=1}^{m}|\alpha_{im}|^2,\nonumber
\end{eqnarray}where $\alpha_{ij}=\alpha(e_{i},e_{j})$.
Let $f:M^{m}\to \mathbb{R}$  given by $f(x)=\vert \varphi
(x)\vert^{2}$. The Hessian  of $f$ at $x \in M$ and $\nu\in T_{x}M$,
$\vert \nu \vert =1$ satisfies
\begin{eqnarray}\label{eqHess2} \frac{1}{2}\Hess f (x)(\nu,\nu)& =&
+\langle \varphi(x),\,\alpha (\nu,\nu)\rangle \nonumber \\
&\geq & \left[1-\vert \varphi (x)\vert \cdot \vert \alpha_{mm}
 \vert \right]\nonumber \\
& \geq & \left[1-\vert \varphi(x)\vert \sqrt{-{Ric (x)(\nu )}} \right] \\
& \geq & \left[1-\rho(x) \sqrt{-{Ric (x)(\nu )}} \right]\nonumber\\
& \geq & 1-c\nonumber
\end{eqnarray} The proof of Theorem (\ref{theorem1})  from equation
(\ref{eqHess}) on  shows that $M$ is is properly immersed in Theorem
(\ref{theorem2}). To show that $M$ has finite topological type
observe that $$ \vert \alpha \vert (\xi (s,y)) \leq \sqrt{-Ric (\xi
(s,y))(\nu, \nu)}\leq c/\rho (\xi(s,y))\leq c/ R(s,y)$$ after
(\ref{eqdf6}) and  we still have (\ref{eqdf7}).

\end{document}